\documentclass{amsproc}
\usepackage[utf8]{inputenc}

\usepackage{amsmath,amssymb,amsthm,mathrsfs,mathtools,accents}
\usepackage{hyperref}
\usepackage{xcolor}

\usepackage{tikz}

\numberwithin{equation}{section}
\newtheorem{theorem}{Theorem}[section]

\theoremstyle{definition}
\newtheorem{definition}[theorem]{Definition}

\title{The GeometricDecomposability package for Macaulay2}

\author[M. Cummings]{Mike Cummings}
\address{Department of Mathematics \& Statistics, McMaster University, Hamilton, ON L8S 4L8, Canada}
\email[Mike Cummings]{cummim5@mcmaster.ca}

\author[A. Van Tuyl]{Adam Van Tuyl}
\address{Department of Mathematics \& Statistics, McMaster University, Hamilton, ON L8S 4L8, Canada}
\email[Adam Van Tuyl]{vantuyl@math.mcmaster.ca}

\date{\today}

\begin{document}

\begin{abstract}  
Using the geometric vertex decomposition property, 
first defined by Knutson, Miller, and Yong, a recursive definition for 
geometrically vertex decomposable ideals was given by
Klein and Rajchgot.
We introduce the 
{\it Macaulay2} package
{\tt GeometricDecomposability} which provides  a suite of tools 
to  experiment and test the geometric vertex decomposability property of an
ideal.
\end{abstract}


\maketitle

%%%%%%%%%%%%%%%%%%%%%%%%%%%%%%%%%%%%%%%%%%%%%%%%%%%%%%%%%%

\section{Introduction}
The {\it geometric vertex decomposition} of an ideal
was first introduced by Knutson, Miller, and Yong
\cite{KMY} as part of their study of vexillary matrix Schubert varieties.  Geometric vertex decomposition
can be viewed as a generalization of a vertex decomposition
of a simplicial complex.  Using the notion of
a geometric vertex decomposition, Klein and Rajchgot
\cite{KR} introduced
{\it geometrically vertex decomposable ideals}.  These
ideals, which are defined recursively,  were partially
inspired by the definition of a vertex decomposable
simplicial complex, a recursively defined family of 
simplicial complexes.

As shown by both \cite{KMY} and \cite{KR}, ideals that
have a geometric vertex decomposition, or are
geometrically vertex decomposable, have other
desirable properties.  As one such example, Klein and Rajchgot 
have shown \cite[Corollary 4.8]{KR} that 
homogeneous geometrically vertex decomposable ideals
are {\it glicci}, i.e., these ideals belong to the Gorenstein liaison class of a complete intersection.   
Further 
properties of  geometrically vertex decomposable ideals
have been developed in  \cite{CDSRVT,DSH,K,KW}.

Due to their recent introduction, there are many features
of geometrically vertex decomposable ideals that are still not
known.  To facilitate further experimentation and exploration,
we have created {\tt GeometricDecomposability}, a package
for {\it Macaulay2} that enables researchers to test and search
for ideals that are geometrically vertex decomposable.  In particular, our package allows the user to check if a given ideal satisfies the geometric vertex decomposition property
of \cite{KMY} or the geometrically vertex decomposable property (or its variants) as found in \cite{KR}. 
Our package can
be found at:
\url{https://macaulay2.com/doc/Macaulay2/share/doc/Macaulay2/GeometricDecomposability/html/index.html}
%\url{https://www.github.com/mgacummings/GeometricDecomposability}.
This note reviews the needed mathematical
background, summarizes
the main features of our packages, and provides
some illustrative examples.

%%%%%%%%%%%%%%%%%%%%%%%%%%%%%%%%%%%%%%%%%%%%%%%%

\section{Mathematical background}
We summarize the mathematical background to define 
the geometric vertex decomposition property and 
geometrically  vertex decomposable ideals. 
Throughout, $k$ denotes a field.

Let $y$ be a variable of the polynomial ring  $R = k[x_1,\ldots,x_n]$. 
A monomial ordering $<$ on $R$ is said to be {\it $y$-compatible} 
if the initial term of $f$ satisfies ${\rm in}_<(f) = 
{\rm in}_<({\rm in}_y(f))$ for all $f \in R$.  Here, 
${\rm in}_y(f)$ is the {\it initial $y$-form} of $f$, that is, 
if $f = \sum_i \alpha_iy^i$ and $\alpha_d \neq 0$ but $\alpha_t = 0$ 
for all $t >d$, then ${\rm in}_y(f) = \alpha_d y^d$.  We set 
${\rm in}_y(I) = \langle {\rm in}_y(f) ~|~ f \in I \rangle$ to 
be the ideal generated by all the initial $y$-forms of an ideal $I$ in $R$.

Given an ideal $I \subseteq R$ and a $y$-compatible monomial ordering $<$, 
let $G(I) = \{ g_1,\ldots,g_m\}$ be a Gr\"obner basis of $I$ 
with respect to this ordering.  For $i=1,\ldots,m$, write $g_i$ as 
$g_i = q_iy^{d_i} + r_i$, where $y$ does not divide any term of $q_i$ and ${\rm in}_y(g_i) = q_iy^{d_i}$.
This second condition is equivalent to no term
of $r_i$ is divisible by $y^{d_i}$.  Given this setup, we define two ideals:
$$C_{y,I} = \langle q_1,\ldots,q_m\rangle
~~\mbox{and}~~N_{y,I} = \langle q_i ~|~ d_i = 0 \rangle.$$
Following Knutson, Miller, and Yong \cite{KMY}, we make the following definition:

\begin{definition}[{\cite[Section 2]{KMY}}]\label{onestep}
With the notation as above, 
the ideal $I$ 
has a {\it geometric
 vertex decomposition with 
 respect to $y$} if
$${\rm in}_y(I) = C_{y,I} \cap (N_{y,I} + \langle y \rangle).$$
\end{definition}

Using Definition \ref{onestep}, Klein and Rajchgot \cite{KR}
recursively defined geometrically vertex decomposable ideals.
Recall that an ideal $I$ is {\it unmixed} if all of its associated 
primes have the same height.

\begin{definition}[{\cite[Definition 2.7]{KR}}]\label{gvd}
An ideal $I$ of $R =k[x_1,\ldots,x_n]$ is {\it geometrically vertex 
decomposable} if $I$ is unmixed and
\begin{enumerate}
    \item[$(1)$]  $I = \langle 1 \rangle$, or $I$ is generated by 
    a (possibly empty) subset of variables of $R$, or
    \item[$(2)$] there is a variable $y = x_i$ in $R$ and a 
    $y$-compatible monomial ordering $<$ such that
    $I$ has a geometric vertex decomposition with respect to $y$, 
    i.e.,
    $${\rm in}_y(I) = C_{y,I} \cap (N_{y,I} + \langle y \rangle),$$
    and the contractions of the ideals $C_{y,I}$ and $N_{y,I}$ to 
    the ring $k[x_1,\ldots,\hat{y},\ldots,x_n]$ are 
    geometrically vertex decomposable.
\end{enumerate}
The ideals $\langle 0 \rangle$ and 
$\langle 1 \rangle$ in the ring $k$ are also considered 
geometrically vertex decomposable by convention.
\end{definition}

Klein and Rajchgot also introduced two variants on the 
geometrically vertex decomposable property.  We also 
recall these definitions.  For the first variant, observe that in
the Definition \ref{gvd},
you do not need to use the induced monomial order for the contractions.  
Indeed, it could be the case that you need to pick different
monomial orders to verify that $C_{y,I}$ and $N_{y,I}$ are geometrically
vertex decomposable.  For $<$-compatibly geometrically vertex decomposable
ideals, a fixed lexicographical order (and its induced monomial
orders) are used throughout; the formal
definition for this class of ideals is given below.

\begin{definition}[{\cite[Definition 2.11]{KR}}]\label{compatible}
Fix a lexicographical order $<$ on $R = k[x_1,\ldots,x_n]$.  
An ideal $I \subseteq R$ is {\it $<$-compatibly geometrically vertex decomposable} if $I$ is unmixed and
\begin{enumerate}
    \item[$(1)$]  $I = \langle 1 \rangle$, or $I$ is generated by 
    a (possibly empty) subset of variables of $R$, or
    \item[$(2)$] for the largest variable $y=x_i$ in $R$ with
    respect to the order $<$, the ideal
    $I$ has a geometric vertex decomposition with respect to $y$, 
    and the contractions of the ideals $C_{y,I}$ and $N_{y,I}$ to 
    the ring $S = k[x_1,\ldots,\hat{y},\ldots,x_n]$ are 
    $<$-compatible geometrically vertex decomposable, where we 
    use $<$ to also denote the  natural induced monomial order on $S$.
\end{enumerate}
\end{definition}

The second variant relaxes some conditions on the ideals
$C_{y,I}$ and $N_{y,I}$, giving a weaker version of 
geometrically vertex decomposable ideals.  The definition
was inspired by Nagel and R\"omer's notion of a
weakly vertex decomposable simplicial complex
(see \cite[Definition 3.1]{NR}).
Following \cite[Section 2]{KR}, the geometric vertex decomposition ${\rm in}_y(I) = C_{y, I} \cap (N_{y, I} + \langle y \rangle)$ is {\em degenerate} if $C_{y, I} = \langle 1 \rangle$ or if $\sqrt{C_{y, I}} = \sqrt{N_{y,I}}$, and it is {\em nondegenerate} otherwise.
Note that the definition requires the field to be infinite.

\begin{definition}[{\cite[Definition 4.6]{KR}}]\label{weak}
Let $k$ be an infinite field.
An ideal $I \subseteq R = k[x_1, \ldots, x_n]$ is {\it weakly geometrically vertex decomposable} if $I$ is unmixed and
\begin{enumerate}
\item[$(1)$] $I = \langle 1 \rangle$, or $I$ is generated by 
a (possibly empty) subset of variables of $R$, or
\item[$(2)$] (Degenerate Case) for some variable 
$y = x_j$ of $R$, 
${\rm in}_y(I) = C_{y,I} \cap (N_{y,I} + \langle y \rangle)$ is
a degenerate geometric vertex decomposition and the contraction 
of $N_{y,I}$ to the ring $k[x_1,\ldots,\hat{y},\ldots,x_n]$
is weakly geometrically vertex decomposable, or
\item[$(3)$] (Nondegenerate Case) for some variable $y = x_j$ of $R$,  ${\rm in}_y(I) = C_{y,I} \cap (N_{y,I} + \langle y \rangle)$ is
a nondegenerate geometric vertex decomposition, the contraction of $C_{y,I}$ to the ring  $k[x_1,\ldots,\hat{y},\ldots,x_n]$
is weakly geometrically vertex decomposable, and $N_{y,I}$ is radical and Cohen-Macaulay.
\end{enumerate}
\end{definition}

Properties of geometrically vertex decomposable ideals,
$<$-compatibly geometrically vertex decomposable ideals,
and weakly geometrically vertex decomposable ideals 
were developed in \cite{KR}.   In particular, it was
shown that these ideals can give insight into questions
related to liaison theory.

%%%%%%%%%%%%%%%%%%%%%%%%%%%%%%%%%%%%%%%%%%%%%%%%%%%%%%%
                        	
\section{The package}

The {\it Macaulay2} package {\tt GeometricDecomposability} was
created as a tool to determine whether an 
ideal $I$ in $R = k[x_1,\ldots,x_n]$ is geometrically
vertex decomposable (or if it satisfies one of the variants). We highlight some of the key features of this package.

Our package is built around the function {\tt oneStepGVD},
which is designed to test whether or not an
ideal $I$ has a geometric vertex decomposition 
with respect to a given variable $y$.  In other words,
this function determines if a given $I$ and $y$ satisfy the properties of Definition \ref{onestep}.  As seen from
Definitions \ref{gvd}, \ref{compatible}, and \ref{weak}, 
checking whether or not an ideal has a geometric
vertex decomposition is key step in these recursive 
definitions.   Our choice of name for this function was 
motivated by the fact that this function
allows us to move one iteration, or ``one step'', in the recursive
definition.

For a given variable $y$, the function {\tt oneStepGVD} 
returns a sequence, where the first element in the sequence 
is true or false depending on whether or not the 
given variable $y$ gives a 
geometric vertex decomposition of $I$, while the second 
element is the ideal $C_{y,I}$, and the third element is the ideal $N_{y,I}$.
As an illustration, we consider the ideal
found in \cite[Example 2.16]{KR}:
\begin{verbatim}
i1 : loadPackage "GeometricDecomposability";
i2 : R = QQ[a..f];
i3 : I = ideal(b*(c*f - a^2), b*d*e, d*e*(c^2+a*c+d*e+f^2));  
i4 : oneStepGVD(I,b)
                              2       2 2        2        2        
o4 = (true, ideal (a*c*d*e + c d*e + d e  + d*e*f , d*e, a  - c*f),
     --------------------------------------------------------------------------
                      2       2 2        2
     ideal(a*c*d*e + c d*e + d e  + d*e*f ))
\end{verbatim}
In this case, we do have a geometric vertex decomposition.
If, on the other hand, we asked if the ideal
has a geometric vertex decomposition with respect to
the variable $c$, we get a negative answer:
\begin{verbatim}
i5 : oneStepGVD(I,c)
o5 = (false, ideal (b*d*e, b*f, a*d*e + d*e), ideal(b*d*e))
\end{verbatim}

We want to highlight that the ideals $C_{y,I}$ and $N_{y,I}$ do not depend upon the choice of the Gr\"obner basis or
a particular $y$-compatible order 
(see the comment after \cite[Definition 2.3]{KR}).
In our package, when we compute $C_{y,I}$ and $N_{y,I}$, we 
use a lexicographical ordering on $R$ where $y > x_j$ for all 
$i \neq j$ if $y = x_i$ since this gives us an easily accessible $y$-compatible 
order.

If the user only requires the ideal
$C_{y,I}$ or $N_{y,I}$, we have built functions to find
these ideals, namely ${\tt CyI}$ and ${\tt NyI}$.  These
functions actually call {\tt oneStepGVD}, and then return
the relevant item in the list.  Continuing with the example
above, we have:
\begin{verbatim}
i6 : CyI(I,c)
o6 = ideal (b*d*e, b*f, a*d*e + d*e)
\end{verbatim}
Note that this is the second entry in the sequence in output {\tt o5}.

As a tool to encourage experimentation, we have
also included the function {\tt findOneStepGVD}.  Given
an ideal $I$ in $R = k[x_1,\ldots,x_n]$, it returns
a list of all the variables with which $I$ has
a geometric vertex decompositon.  In our running example,
there is only one such variable:
\begin{verbatim}
i7 : findOneStepGVD(I)
o7 = {b}
\end{verbatim}

We have created three separate functions
to check if an ideal is geometrically
vertex decomposable, $<$-compatibly 
geometrically vertex decomposable, or weakly
geometrically vertex decomposable.  
The implementation of each function
requires repeated use of the function
{\tt oneStepGVD}.

To show that our running
example is geometrically vertex
decomposable, we enter the 
command:
\begin{verbatim}
i8 : isGVD(I)
o8 = true
\end{verbatim}
Running the above command with 
the 
option {\tt isGVD(I,Verbose=>True)} will allow the user
to identify which  variable is used
to form the geometric vertex decomposition at
each step. 

Using the
{\tt isLexCompatiblyGVD} command,
we can check if an ideal is
$<$-compatibly geometrically vertex 
decomposable.  The user must specify the ideal,
and the specific
lexicographical order by providing
an ordering of the variables.  As an
example, we can check if our running
example is $<$-compatibly geometrically vertex
decomposable with respect to the 
lexicographical order which satisfies
$f>e>d>c>b>a$.  The specific command is:
\begin{verbatim}
i9 : isLexCompatiblyGVD(I,{f,e,d,c,b,a})
o9 = false
\end{verbatim}
We can also search over all possible lexicographical
orders to determine if an ideal is $<$-compatibly
geometrically vertex decomposable.  Specifically, the command {\tt findLexCompatiblyGVDOrders}
returns all the lexicographical orders for which the
ideal is $<$-compatibly geometrically vertex decomposable.
In our example, there is no
such lexicographical order, as given by the output:
\begin{verbatim}
i10 : findLexCompatiblyGVDOrders(I)
o10 = {}
\end{verbatim}
This agrees with \cite[Exampe 2.16]{KR} which proved that 
this ideal is not $<$-compatibly geometrically vertex 
decomposable.  Note that running this command 
can be  computationally expensive since 
one may need to check $n!$ different 
lexicographical orders in $k[x_1,\ldots,x_n]$.

Finally, the command {\tt isWeaklyGVD} enables
the user to check if an ideal satisfies
Definition \ref{weak}.  By \cite[Corollary 4.7]{KR}, 
all geometrically vertex decomposable ideals are weakly
geometrically vertex decomposable.   Our running 
example is therefore weakly geometrically vertex decomposable,
as expected:
\begin{verbatim}
i11 : isWeaklyGVD(I)
o11 = true
\end{verbatim}
There are ideals that are weakly geometricaly
vertex decomposable, but not geometrically vertex
decomposable.  The following example comes
from \cite[Example 4.10]{KR}:
\begin{verbatim}
i12 : J = ideal(b*(c*f-a^2),b*d*e,d*e*(a^2+f^2+d*e));
i13 : isWeaklyGVD(J)
o13 = true
i14 : isGVD(J)
o14 = false
\end{verbatim}
We present a new example of an ideal with this property 
at the end of this paper.
\section{Examples}
We illustrate the {\tt GeometricDecomposability} package using examples 
from square-free monomial ideals and toric ideals of graphs. 
\subsection{Square-free monomial ideals} As noted in the introduction,
geometrically vertex decomposable ideals was inspired by
the definition of vertex decomposable simplicial complexes.  
We flesh out this connection, and explain the output
of package within the context of square-free monomial ideals.

Let $V = \{x_1,\ldots,x_n\}$ be a vertex set, and let $2^V$ denote
the power set of $V$.  A {\it simplicial complex} is a 
subset $\Delta \subseteq 2^V$ that satisfies the two properties:
$(i)$ if $F \in \Delta$ and $G \subseteq F$, then $G \in \Delta$, and
$(ii)$ $\{x_i\} \in \Delta$ for all $i \in \{1,\ldots,n\}$.  The 
maximal elements of $\Delta$ with respect to inclusion are called the
{\it facets}.  A simplicial complex is {\it pure} if all the facets
have the same cardinalty.  Given a vertex $x \in V$, the 
{\it deletion} of $x$ is the subcomplex ${\rm del}_\Delta(x) = \{F \in 
\Delta ~|~ F \cap \{x\} = \emptyset\}$ and the {\it link} of
$x$ is ${\rm link}_\Delta(x) = \{F ~|~ F \cap \{x\} = \emptyset
~~\mbox{and}~~ F \cup \{x\} \in \Delta\}$.
Note that the link and
deletion are not necessarily simplicial complexes on $V$, but on a subset of $V$.  
Precisely, ${\rm link}_\Delta(x)$ is 
simplicial complex on $\bigcup_{F \in {\rm link}_\Delta(x)} F$ and ${\rm del}_\Delta(x)$
is a simplicial complex on 
$\bigcup_{F \in {\rm del}_\Delta(x)} F$.

Vertex decomposable simplicial complexes were first
introduced  by Provan and Billera \cite{PB}:
\begin{definition}\label{vd}
A simplicial complex $\Delta$ on $V$ is {\it vertex decomposable}
if $\Delta$ is pure and either $(i)$ $\Delta = \emptyset$, or 
$\Delta$ is a simplex, i.e., the only facet is
$\{x_1,\ldots,x_n\}$, or $(ii)$, there exists a vertex
$x \in V$ such that ${\rm del}_\Delta(x)$ and ${\rm link}_\Delta(x)$ are vertex decomposable.
\end{definition}
\noindent

To connect Definition \ref{vd} with the
definition of geometrically vertex decomposable ideals
in Definition \ref{gvd}, we use the Stanley-Reisner
correspondence.  In particular, given 
a simpicial complex $\Delta$, define
$$I_\Delta = \langle x_{i_1}\cdots x_{i_j} ~|~
\{x_{i_1},\ldots,x_{i_j}\} \not\in \Delta \rangle
\subseteq R = k[x_1,\ldots,x_n]$$
to be the square-free monomial ideal generated by monomials
corresponding to subsets not in $\Delta$.  The connection
between the two definitions then comes via the following
theorem:

\begin{theorem}[{\cite[Proposition 2.9]{KR}}]\label{vd=gvd}Let
$\Delta$ be a simplicial complex on 
$V = \{x_1,\ldots,x_n\}$.  Then $\Delta$ is 
vertex decomposable if and only if $I_\Delta$
is geometrically vertex decomposable.
\end{theorem}
\noindent
In other words, 
the above theorem says that the square-free
monomial ideals that are geometrically vertex decomposable
are precisely those square-free monomial ideals that
are the Stanley-Reisner ideals of a vertex decomposable
simplicial complex.

Our example below highlights the connection
between the  ${\rm link}_\Delta(x)$ and
${\rm del}_\Delta(x)$ and the ideals $C_{y,I}$
and $N_{y,I}$ that appear in the 
definition of (weakly) 
geometrically vertex decomposable ideals.   
Starting a fresh \emph{Macaulay2} session, consider
the simplicial complex $\Delta$ on the vertex
set  $\{a,b,c,d,e\}$ with facets 
$\{\{a,c\},\{a,d\},\{b,d\},\{b,e\},\{c,e\}\}$.  
Using the
{\tt SimplicialDecomposability} package \cite{CII}, we can input this simplicial
complex and check that it is 
vertex decomposable:
\begin{verbatim}
i1 : loadPackage "GeometricDecomposability";
i2 : loadPackage "SimplicialDecomposability";
i3 : R = QQ[a..e];
i4 : Delta = simplicialComplex {a*c,a*d,b*d,b*e,c*e};
i5 : isVertexDecomposable(Delta)
o5 = true
\end{verbatim}
We can now obtain the simplicial complexes
${\rm del}_\Delta(a)$ and ${\rm link}_\Delta(a)$ 
of the vertex $a \in \{a,\ldots,e\}$ using the
following commands.  We also include the corresponding Stanley-Reisner 
ideal of each simplicial complex:
\begin{verbatim}
i6 : Link = link(Delta,a);
i7 : Delete = faceDelete(a,Delta);
i8 : IDelta = monomialIdeal Delta
o8 = monomialIdeal (a*b, b*c, c*d, d*e, e*a)
i9 : monomialIdeal Link
o9 = monomialIdeal (a, b, c*d, e)
i10 : monomialIdeal Delete
o10 = monomialIdeal (a, b*c, c*d, d*e)
\end{verbatim}
Now consider the output of the {\tt oneStepGVD}
command with input $I_\Delta$ and the vertex $a$:
\begin{verbatim}
i11 : oneStepGVD(IDelta,a)
o11 = (true, ideal (d*e, c*d, b*c, e, b), ideal (d*e, c*d, b*c))
\end{verbatim}
If we compare the ideals $C_{a,I_\Delta}$
and $N_{a,I_\Delta}$, the second and third ideals in 
the above list, with the Stanley-Reisner ideals
of ${\rm link}_\Delta(a)$ and ${\rm del}_\Delta(a)$, they
are the same except that the later ideals have an extra 
generator, namely the variable $a$. (Note that the
given generators of $C_{a,I_\Delta}$ are not a minimal set of generators).  Technically, the simplicial complexes
${\rm link}_\Delta(a)$ and ${\rm del}_\Delta(a)$ are 
simplicial complexes on the vertex set $\{b,c,d,e\}$,
and so their corresponding Stanley-Reisner ideals
should belong to $k[b,c,d,e]$.  We can move all the
ideals to this ring, and then we verify that 
we have an equality of ideals:
\begin{verbatim}
i12 : S = QQ[b,c,d,e]
i13 : C1=substitute(CyI(IDelta,a),S)
o13 = ideal (d*e, c*d, b*c, e, b)
i14 : N1=substitute(NyI(IDelta,a),S)
o14 = ideal (d*e, c*d, b*c)
i15 : L1=substitute(monomialIdeal Link,S)
o15 = ideal (0, b, c*d, e)
i16 : D1=substitute(monomialIdeal Delete,S)
o16 = ideal (0, b*c, c*d, d*e)
\end{verbatim}
In general, the ideal $C_{y,I}$, respectively the ideal
$N_{y,I}$, can be viewed as the algebraic analog of the link of
a vertex, respectively the deletion of a vertex, of a simplicial complex.

As a final calculation, we verify that our Stanley-Reisner ideal 
is geometrically vertex decomposable, as
expected by Theorem \ref{vd=gvd}:
\begin{verbatim}
i17 : isGVD(IDelta)
o17 = true
\end{verbatim}

\subsection{Toric ideals of graphs}
For our second example, we use our package to find minimal examples
of toric ideals of graphs that are (weakly) geometrically vertex decomposable.

Let $G = (V,E)$ be a finite simple graph with
vertex set $V = \{x_1,\ldots,x_n\}$ and edge set $E = \{e_1,\ldots,e_m\}$.
Let $R = k[e_1,\ldots,e_m]$ and $S = k[x_1,\ldots,x_n]$.  We define a map
$\varphi:R \rightarrow S$ by $\varphi(e_i) = x_jx_k$ where $e_i = \{x_j,x_k\}$. 
The \emph{toric ideal of $G$}, denoted $I_G$, is the kernel of this map,
that is, $I_G = {\rm ker}\varphi$.   It can be shown 
(see, for example \cite[Proposition 10.1.5]{V}) that $I_G$ is prime
binomial ideal.  Furthermore, the generators of $I_G$ correspond to 
closed even walks in the graph $G$.  Informally, 
a closed even walk is a sequence
of adjacent edges that start and stop at the same edge (see
\cite[Chapter 7.1]{V} for the formal definition).
Not every graph has a closed even walk, e.g., trees, so in some cases $I_G = \langle 0 \rangle$.

Geometrically vertex decomposable toric ideals of graphs were first
studied in \cite{CDSRVT}.  It was shown that every
toric ideal of a bipartite graph is  geometrically
vertex decomposable.  However, not every toric ideal of graph
is geometrically vertex decomposable, as noted 
in \cite[Remark 7.1]{CDSRVT}.

Using our package {\tt GeometricDecomposability} we can find 
minimal examples of graphs that are (weakly) geometric vertex 
decomposable.  By \cite[Theorem 3.3]{CDSRVT}, we can restrict our search to connected graphs.
Table \ref{table summary} summarizes the results
of our computation.  For all connected graphs with $e$
edges with $e \in \{4,\ldots,9\}$, we first constructed the toric
ideal $I_G$.  If $I_G \neq \langle 0 \rangle$, we then checked
if the ideal was geometrically vertex decomposable and weakly
geometrically vertex decomposable.  The second column of 
Table \ref{table summary} is the number of simple connected
graphs on $e$ edges (this is sequence
A002905 in \cite{OEIS}), the third column is the number of such
graphs that have a non-zero toric ideal, while the fourth
and fifth record the number of these ideals that 
are geometrically vertex decomposable (GVD) or weakly GVD.

This data was computed also using the \emph{Macaulay2} packages {\tt Nauty} and {\tt FourTiTwo} \cite{Nauty, FourTiTwo}.
We used the following code:

\footnotesize
\begin{verbatim}
loadPackage "GeometricDecomposability";
loadPackage "NautyGraphs";
loadPackage "FourTiTwo";

getToricIdeal = (A,R) -> (
    -- A, an incidence matrix; R, a ring
    m = product gens R;
    return saturate(sub(toBinomial(transpose(syz(A)),R),R),m);
)

R = QQ[a..i]; -- ring with 9 indeterminates (for the up to 9 edges we will see)

-- all graphs with number of edges between 4 and 9 (inclusive)
GList = flatten for n from 4 to 10 list generateGraphs(n, 4, 9, OnlyConnected=>true);

-- for each number of edges, filter the list to these graphs, and check if GVD
for E from 4 to 9 do (
    f := buildGraphFilter {"NumEdges" => E};
    HList := filterGraphs(GList, f);
    print("There are " | toString (#HList) | " graphs on " | toString E | " edges,");

    -- create lists of graphs, incidence matrices of these graphs, and toric ideals of these graphs
    graphList = for g in HList list stringToGraph g;
    incMatList = for G in graphList list incidenceMatrix G;
    idealList = for A in incMatList list getToricIdeal A;

    -- now look at the subset of the ideals which are not the zero idea, weakly GVD, and GVD
    nonZeroList = select(idealList, i-> i != 0);
    print(toString(#nonZeroList) | " of which are nonzero");

    wgvdList = select(nonZeroList, i -> isWeaklyGVD(i));
    print(toString(#wgvdList) | " of which are weakly GVD and");

    gvdList = select(nonZeroList, i -> isGVD(i));
    print(toString(#gvdList) | " of which are GVD");
    print("");
);
\end{verbatim}
\normalsize
\renewcommand{\arraystretch}{1.05}
\begin{table}
\centering
    \begin{tabular}{p{1.1cm}p{2.75cm}p{2.15cm}p{2.05cm}p{2.95cm}} 
        \hline
        \hline 
        edges & simple connected graphs & 
        non-zero toric ideals &
        weakly GVD toric ideals &
        GVD toric ideals
        \\ \hline 
        \hline
        4 & 5 & 1 & 1 & 1 \\ %\hline 
        5 & 12 & 2 & 2 &2  \\ %\hline 
        6 & 30 & 11 & 11 &  11\\ %\hline 
        7 & 79 & 33 & 33 & 33  \\ %\hline 
        8 & 227 & 125 & 124 & 124\\ %\hline
        9  & 710 & 449 & 445 & 444 \\
        \hline
        \hline
    \end{tabular}
    \vspace{.2cm}
    \caption{Comparison of the  number of (weakly) GVD toric
    ideals of graphs}
    \label{table summary}
\end{table}

\begin{figure}[!ht]
    \centering
    \begin{tikzpicture}[scale=0.27]
      % graph G
      
      \draw (-14,3) -- (-14,9) node[midway,  left] {$a$};
      \draw (-10.5,6) -- (-14,9) node[midway, above] {$b$};
      \draw (-14,3) -- (-10.5,6) node[midway, below] {$c$};
      \draw (-10.5,6) -- (-7.5,6) node[midway, above] {$d$};
      \draw (-7.5,6) -- (-4.5,6) node[midway, above] {$e$};
      \draw (-4.5,6) -- (-1,3) node[midway, below] {$f$};
      \draw (-4.5,6) -- (-1,9) node[midway, above] {$g$};
      \draw (-1,3) -- (-1,9) node[midway, right] {$h$};

      \fill[fill=white,draw=black] (-14,3) circle (.2) node[left]{};
      \fill[fill=white,draw=black] (-14,9) circle (.2) node[left]{};
      \fill[fill=white,draw=black] (-10.5,6) circle (.2) node[below]{};
      \fill[fill=white,draw=black] (-7.5,6) circle (.2) node[below]{};
      \fill[fill=white,draw=black] (-4.5,6) circle (.2) node[below]{};
      \fill[fill=white,draw=black] (-1,3) circle (.2) node[right]{};
      \fill[fill=white,draw=black] (-1,9) circle (.2) node[right]{};

      \draw (6,3) -- (6,9) node[midway,  left] {$a$};
      \draw (9.5,6) -- (6,9) node[midway, above] {$b$};
      \draw (6,3) -- (9.5,6) node[midway, below] {$c$};
      \draw (9.5,6) -- (12.5,6) node[midway, above] {$d$};
      \draw (12.5,6) -- (15.5,3) node[midway, below] {$f$};
      \draw (12.5,6) -- (15.5,9) node[midway, above] {$e$};
      \draw (15.5,3) -- (15.5,9) node[midway, right] {$g$};
      \draw (15.5,9) -- (18.5,6) node[midway, right] {$h$};
      \draw (15.5,3) -- (18.5,6) node[midway, right] {$i$};

      \fill[fill=white,draw=black] (6,3) circle (.2) node[left] {};
      \fill[fill=white,draw=black] (6,9) circle (.2) node[left] {};
      \fill[fill=white,draw=black] (9.5,6) circle (.2) node[below] {};
      \fill[fill=white,draw=black] (12.5,6) circle (.2) node[below] {};
      \fill[fill=white,draw=black] (15.5,3) circle (.2) node[right] {};
      \fill[fill=white,draw=black] (15.5,9) circle (.2) node[right] {};
      \fill[fill=white,draw=black] (18.5,6) circle (.2) node[right] {};
     
    \end{tikzpicture}
       \caption{The unique graph $G$ with 8 edges whose toric ideal is not GVD (on the left) and the unique
       graph $H$ with 9 edges whose toric ideal is weakly GVD but not GVD (on the right)}
    \label{fig:2triangles}
\end{figure}
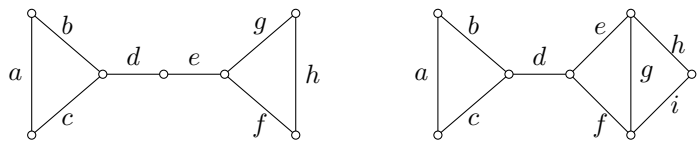

Table \ref{table summary} implies that there is exactly one
graph on 8 edges (and hence the smallest graph) whose toric ideal is not geometrically vertex decomposable.  
Figure \ref{fig:2triangles} shows the unique graph $G$ on 8 edges whose toric ideal is not geometrically vertex decomposable.  
Note that this graph is
the same graph of \cite[Remark 7.1]{CDSRVT}.
The toric ideal of this graph is $I_G = \langle ad^2fg - bce^2h \rangle$.

Our computations also imply there is a unique
graph on 9 edges whose toric ideal is weakly
geometrically vertex decomposable, but not geometrically
vertex decomposable (and moreover, this is the 
smallest such graph). This graph $H$ appears to the
right in Figure \ref{fig:2triangles}.
The toric ideal of this graph is $I_H = \langle fg-ei, bcef - ad^2g \rangle$.

\noindent
{\bf Acknowledgments.} 
We thank Sergio Da Silva, Patrica Klein, and Jenna Rajchgot
for comments and suggestions.
We additionally thank the anonymous referees for concrete suggestions that improved both this manuscript and the software package.
Cummings was 
partially supported by the Milos Novotny Fellowship and an NSERC USRA and CGS-M.  
Van Tuyl’s research is partially
supported by NSERC Discovery Grant 2019-05412. 

\bibliographystyle{plain}
\bibliography{references}

\end{document}